\documentclass[11pt]{imsart}
\usepackage{amsmath,amssymb,amsthm}
\usepackage[utf8]{inputenc}

\usepackage[english]{babel}

\usepackage{mathrsfs, euscript}
\usepackage{graphicx}
\usepackage{stackrel}
\usepackage{gensymb}

\usepackage[round]{natbib}

\begin{document}

\renewcommand{\baselinestretch}{1.2}

\newcommand{\cc}{\c{c}}
\newcommand{\ic}{\^{\i}}
\newcommand{\al}{\alpha}
\newcommand{\om}{\omega}
\newcommand{\Om}{\Omega}
\newcommand{\La}{\Lambda}
\newcommand{\la}{\lambda}
\newcommand{\De}{\Delta}
\newcommand{\de}{\delta}
\newcommand{\ep}{\epsilon}
\newcommand{\be}{\beta}
\newcommand{\ga}{\gamma}
\newcommand{\Ga}{\Gamma}
\newcommand{\te}{\theta}
\newcommand{\Te}{\Theta}
\newcommand{\si}{\sigma}
\newcommand{\Si}{\Sigma}
\newcommand{\vep}{\varepsilon}
\newcommand{\veps}{\varepsilon}
\newcommand{\eps}{\epsilon}
\newcommand{\ze}{\zeta}
\newcommand{\T}{\mathrm{\tiny T}} 
\newcommand{\pa}{\partial}
\newcommand{\dd}[2]{\frac{\pa #1}{\pa #2}} 

\newcommand{\Z}{\mathbb{Z}}
\newcommand{\R}{\mathbb{R}}
\newcommand{\CC}{\mathbb{C}}   
\newcommand{\Q}{\mathbb{Q}}
\newcommand{\N}{\mathbb{N}}

\newcommand{\cov}{{\mbox{ Cov} }}
\newcommand{\var}{{\mbox{ Var} }}

\newcommand{\PP}{\mathbb{P}}
\newcommand{\E}{\mathbb{E}}

\newcommand{\rw}{\rightarrow}
\newcommand{\lrw}{\longrightarrow}
\newcommand{\Rw}{\Rightarrow}
\newcommand{\Lrw}{\Longrightarrow}
\newcommand{\bigo}{\bigotimes}
\newcommand{\ti}{\times}
\newcommand{\oti}{\otimes}

\def\b#1{{#1}}
\def\C#1{{\cal #1}} 

\newcommand{\CF}{{\mathscr F}}
\newcommand{\CH}{{\mathscr H}}
\newcommand{\CT}{{\mathscr T}}
\newcommand{\CN}{{\mathscr N}}
\newcommand{\CX}{{\mathscr X}}
\newcommand{\CY}{{\mathscr Y}}
\newcommand{\CZ}{{\mathscr Z}}
\newcommand{\CD}{{\mathscr D}}
\newcommand{\CB}{{\mathscr B}}
\newcommand{\CA}{{\mathscr A}}
\newcommand{\CG}{{\mathscr G}}
\newcommand{\CQ}{{\mathscr Q}}
\newcommand{\CL}{{\mathscr L}}
\newcommand{\CV}{{\mathscr V}}
\newcommand{\CE}{{\mathscr E}}
\newcommand{\CP}{{\mathscr P}}
\newcommand{\CU}{{\mathscr U}}
\newcommand{\bl}{\vspace{1ex}}
\newtheorem{etheo}{Theorem}
\newtheorem{elemm}{Lemma}
\newtheorem{eprop}{Proposition}
\newtheorem{ecoro}{Corollary}
\def\minmax{\mathop{\hbox{\rm min\,max}}}
\def\tr{\mathop{\hbox{\rm tr}}}
\def\Arg{\mathop{\hbox{\rm Arg}}}
\def\diag{\mathop{\hbox{\rm diag}}}
\def\grad{\mathop{\nabla}}
\def\card{\mathop{\hbox{\rm card}}}
\def\bc{\mathop{\hbox{\mbox{\large\rm $|$}}}}
\def\telque{\; : \;} 
\def\pds#1#2{{\langle #1,#2 \rangle}}         
\def\bpv{\mathop{;}}                          
\def\dfrac#1#2{\frac{\displaystyle #1}{\displaystyle #2}}
\def\ffrac#1{\frac{1}{#1}}
\newcommand{\norm}[2]{ \| #1 \|_{#2} }
\def\un{\mbox{\large 1\hskip-0.30em I}} 
\def\indic#1{\un_{#1}}                        
\def\cf{{\em cf. }} 
\def\sgn{\mathop{\hbox{\rm sgn}}}
\newcommand{\rmx}[1]{{\mbox{\rm #1}}}
\newcommand{\rmxs}[1]{{\mbox{\rm {\scriptsize #1}}}}
\newcommand{\rmxt}[1]{{\mbox{\rm {\tiny #1}}}}
\newcommand{\trsp}{ {}^t \! }                
\newcommand{\fq}[2]{\trsp #1 \; #2 \; #1}

\def\demof{\noindent{\textbf{Preuve. \quad}}}
\def\demo{\noindent{\textbf{Proof. \quad}}}
\def\ddemof#1{\noindent{\textbf{Preuve {#1}. \quad}}}
\def\ddemo#1{\noindent{\textbf{Proof {#1}. \quad}}}

\def\findemo{\hskip2mm\mbox{\vbox{\hrule height 3pt depth 6pt width 6pt}}}

%

\newcommand{\eq}[1]{(\ref{#1})}         
\newcommand{\eqn}[1]{Eq.(\ref{#1})}         
\newcommand{\cte}{\rmx{cte }}

\def\remise{\setcounter{equation}{0}}

\def\slabel#1{\label{#1} \hfill \hskip2mm\hbox to 4mm{\dotfill}%
        \mbox{\tiny\em(#1)}\hbox to 4mm{\dotfill}}

\newcounter{exo}
\renewcommand\theexo{\arabic{exo}}
\newenvironment{exo}{\par\refstepcounter{exo}%
\medskip\noindent\textbf{Exercice \theexo}\space}

\newcounter{exoc}
\renewcommand\theexoc{\arabic{chapter}.\arabic{exoc}}
\newenvironment{exoc}{\par\refstepcounter{exoc}%
\medskip\noindent\textbf{\theexoc}\space}

\newcommand\nbin[2]{\begin{pmatrix}#1 \\ #2 \end{pmatrix}}

\newcommand\eurosym{\textgreek{\euro}}

\newtheorem{theorem}{Theorem}
\newtheorem{proposition}{Proposition}
\newtheorem{lemma}{Lemma}
\newtheorem{remark}{Remark}

\renewcommand{\topfraction}{0.85}
\renewcommand{\textfraction}{0.1}
\renewcommand{\floatpagefraction}{0.75}

\newcommand\bveps{{\boldsymbol{\varepsilon}}}
\newcommand\bX{\mathbf{X}}
\newcommand\bY{\mathbf{Y}}
\newcommand\bV{\mathbf{V}}
\newcommand\bD{\mathbf{D}}
\newcommand\bS{\mathbf{S}}
\newcommand\bZ{\mathbf{Z}}
\newcommand\bx{\mathbf{x}}                                            
\newcommand\by{\mathbf{y}}                                            
\newcommand\bz{\mathbf{z}}
\newcommand{\MP}{Mar\v{c}enko-Pastur}

\begin{frontmatter}
  \title{A note on a \MP\ type theorem for time series}
  \runtitle{\MP\ theorem for time series}
  \thankstext{T1}{The research of the author was supported partly by 
    a HKU Start-up grant.}

  \begin{aug}
    \author{\fnms{Jianfeng} \snm{Yao}\thanksref{T1}\ead[label=e3]{jeffyao@hku.hk}}

    \runauthor{J.  Yao}

    \address{Jianfeng Yao \\
      Department of Statistics and Actuarial Science\\
      The University of Hong Kong\\
      Pokfulam, HONG KONG\\
      \printead{e3}
    }
  \end{aug}

  \begin{abstract}
    In this note we develop an extension of the \MP\ 
    theorem to time series model with temporal correlations.
    The limiting spectral distribution (LSD) of the sample covariance
    matrix is characterised by an explicit equation for its Stieltjes
    transform depending on the spectral density of the time series.
    A numerical algorithm is then given to compute the density
    functions of these LSD's. 
  \end{abstract}


\end{frontmatter}

\section{Introduction}

Let $\{\bX_j\}, j=1,\ldots,n$ be  a sequence of $p$-dimensional 
real-valued random vectors and consider the associated 
empirical covariance matrix   
\begin{equation}\label{Sn}
  S_n = \frac1n \sum_{j=1}^n \bX_j \bX_j^\top.  
\end{equation}
The study of the empirical spectral
distribution (ESD) $F_n$    of $S_n$, i.e.  
the distribution generated by its (real-valued) 
eigenvalues, goes back to Wishart in 1920's.  
A milestone work by \citet{MP}
states that if both sample size $n$ and data dimension $p$ 
proportionally grow to infinity  such that $\lim p/n = c$ for some positive $c>0$
and all the coordinates of all the vectors $\bX_j$'s are i.i.d. 
 with mean zero and variance 1, then $F_n$ converges to 
a   nonrandom  distribution. This limiting spectral distribution (LSD),
named after them as the Mar\v{c}enko-Pastur distribution of index $c$ 
has a density function 
\begin{equation}\label{mpdensity}
  f(x)=\frac{1}{2\pi cx} \sqrt{(b-x)(x-a)}, \quad a\le x\le b, 
\end{equation}
with  $a= (1-\sqrt{c})^2$ and
$b=(1+\sqrt{c})^2$ defining the support interval and  has a point mass
$1-1/c$ at the origin if $c>1$. 
Further refinements are made successively by many researchers 
including \citet{Jonsson82},
\citet{Wachter78} and \citet{Yin86}. 

An important work by \citet{Silverstein95} aimed at relaxing  the
independence structure 
between the coordinates of the $\bX_j$'s and considered random vectors
of form $\bY_j=T_p^{1/2}\bX_j$ where $(T_p)$ is a sequence of
non-negative definite matrices.
Assuming that $(T_p)$  is bounded 
in spectral norm and the sequence of ESD of $(T_p)$ has a weak limit 
$H$, he established  a (strong) LSD  for the sample covariance 
matrix $n^{-1}\sum Y_j Y_j^\top$ and provides a characteristic 
equation for  its Stieltjes transform. 
Despite a big step made by this generalisation,  it still does not cover 
all possible correlation patterns of coordinates.
Pursuing these efforts,  a recent work by \citet{BaiZhou08}
pushes  a step further Silverstein's result 
by allowing a very general pattern for correlations between the
coordinates of the $\bX_j$'s satisfying a mild moment conditions. 

In this work, we  extend  such \MP\ type theorems along another
direction by considering time series observations 
instead of an i.i.d. sample. Let us 
first consider an  univariate real-valued linear process 
\begin{equation}\label{zt}
  z_t = \sum_{k=0}^\infty \phi_k \veps_{t-k},\quad t\in\Z,
\end{equation}
where $(\veps_k)$ is a real-valued and weakly stationary white 
 noise with mean zero and variance 1.  The $p$-dimensional process $(\bX_t)$ considered
in this paper will be made by 
$p$ independent copies of the  linear process $(z_t)$, i.e. 
for $\bX_t=(X_{1t},\ldots,X_{pt})^\top$, 
\begin{equation}\label{model}
     X_{it}    = \sum_{k=0}^\infty \phi_k \veps_{i,t-k},\quad t\in\Z,
\end{equation}
where the $p$ coordinate  processes 
$\{(\veps_{1,t},\ldots,\veps_{p,t})\}$  are independent
copies of the univariate error process  $\{\veps_{t}\}$ in \eqref{zt}.
Let $\bX_1,\ldots,\bX_n$ be the observations of  the time series at
time epochs $t=1,\ldots,n$. Again we are interested in the ESD of the
sample covariance matrix $S_n$ in \eqref{Sn}.

The author should mention that a similar problem has been considered 
in  \citet{JWMH09}. However we propose much more general results 
in this note since firstly 
their results are limited to ARMA-type processes instead of a general linear process 
considered here and secondly, they do not find a general equation as
the one proposed in Theorem~\ref{th:main} below except for two 
simplest particular cases of AR(1) and MA(1).

\section{A \MP\ type theorem for linear processes}

Recall that the Stieltjes transform $s_\mu$ of a probability measure
$\mu$ on the real line is a map from the set $\CC^+$ of complex
numbers with positive imaginary part onto itself and defined by 
\[  s_\mu(z)= \int \frac1{x-z}\mu(dx),\quad z\in\CC^+.
\]
We always employ an usual convention that for any complex number $z$, 
$\sqrt{z}$ denotes  its  square
root with nonnegative imaginary part. 

\begin{theorem}\label{th:main}
  Assume that the following conditions hold:
  \begin{enumerate}
  \item The dimensions $p\to\infty$, $n\to\infty$ and $p/n\to c\in (0,\infty)$;
  \item The error process has a fourth moment: $\E \veps_t^4<\infty$; 
  \item The linear filter $(\phi_k)$ is  absolutely
    summable, i.e. $\displaystyle \sum\limits_{k=0}^\infty|\phi_k|<\infty$.
  \end{enumerate}
  Then almost surely  the ESD  of $S_n$ tends to a non-random
  probability distribution $F$. Moreover, the  Stieltjes transform
  $s=s(z)$ of $F$  (as a mapping from $\CC^+$ into  $\CC^+$) 
  satisfies the equation
  \begin{equation}\label{sol}
    z = -\frac1{s}   + \frac{1}{2\pi}\int_{0}^{2\pi}  \frac1{cs+\{2\pi f(\la)\}^{-1}}\,d\la~,
  \end{equation}
  where $f(\la)$ is the  spectral density of the linear process $(z_t)$:
  \begin{equation}\label{sd}
  f(\la)=\frac1{2\pi} \left| \sum_{k=0}^\infty  \phi_k e^{ik\la}   \right|^2,\quad
  \la \in [0,{2\pi}).
  \end{equation}
\end{theorem}

The proof of the theorem is postponed to Section~\ref{proofs}. Let us
mention that although the case $c=0$ is beyond the scope of
Theorem~\ref{th:main}, Equation~\eqref{sol} leads in this case to the 
solution $s(z)=1/(\ga_0-z)$, that is the LSD would be the Dirac mass
at $\ga_0=\var(X_{it})$. This conjectures an extension of
Theorem~\ref{th:main} to the so-called  ``very large $p$ and small
$n$'' asymptotics where one assumes $p\rw\infty$,  $n\rw\infty$ and
$p/n\rw 0$. Indeed, in this scenario taking into account that the
population covariance matrix of $\bX_t$ equals $\ga_0 I_p$, one can
expect that the sample eigenvalues of $S_n$ stay close to the population
ones (all equal to $\ga_0$). Note that such results exist for 
i.i.d. sequence $(\bX_t)$ with i.i.d. components \citep[see][]{BaiYin88b}.

\subsection{Support of the LSD $F$}

Starting from Eq.~\ref{sol} and following the techniques devised in
\citet{SilversteinChoi95}, we can describe precisely the support of the LSD
$F$ in previous theorem.  

\begin{figure}[t]
  \begin{center}
    \includegraphics[width=15pc]{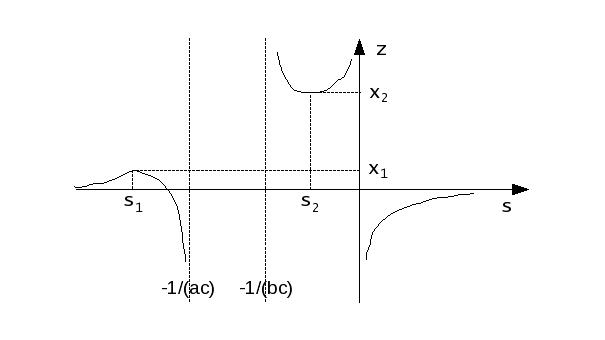}\hskip-1cm
    \includegraphics[width=15pc]{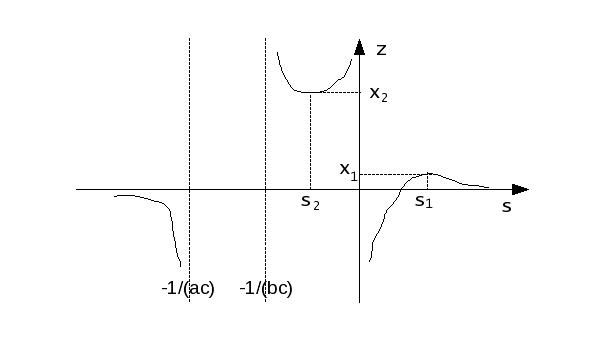}
  \end{center}
  \caption{Determination of the support $[x_1,x_2]$ of the LSD. Left
    panel $c>1$ and right panel $c<1$.\label{fig:map-g}}
\end{figure}

Let $a$ and $b$ be respectively the minimum
and the maximum of the function  $2\pi f$ over $[0,2\pi]$. As $f$ is
infinitely differentiable and positive everywhere,  
both $a$ and $b$ are attained and the range of $2\pi f$ is
exactly the interval $[a,b]$. We will always exclude the situation
$a=0$ since it corresponds to a special class of linear processes, 
namely non-invertible ARMA processes, see \citet[chap.
  9]{GreSze58}, which has no practical interest for applications.  
Therefore the map $s\mapsto z=g(s)$ in Eq.\eqref{sol} has a trace for {\em real-valued} $s$ 
providing $s\notin [-\frac1{ac}, -\frac1{bc}]$. 
Figure~\ref{fig:map-g} depicts this map for both $c<1$ and $c>1$ cases.

The following proposition is a straightforward application of  
results from \citet{SilversteinChoi95} and we then omit its proof.
 
\begin{proposition}\label{prop}
With the map $g:~s\mapsto z=g(s)$ in Eq.\eqref{sol} restricted to real
$s\notin [-\frac1{ac}, -\frac1{bc}]$ (Figure~\ref{fig:map-g})  the following holds:
\begin{enumerate}
\item The LSD $F$  has a compact support $[x_1, x_2]\subset[0,\infty)$
  on which it has a continuous density function.  
  In case of  $c>1$, $F$ has an additional point mass $1-1/c$ at the origin. 
\item When $c>1$,  the map $g$ has an unique maximum $s_1$ on $(-\infty,-\frac1{ac} )$
  and an unique minimum $s_2$ on $(-\frac1{bc},0)$ and we have  
  The edges of the support interval are given by these local extrema:
  $x_1=g(s_1)$ and $x_2=g(s_2)$.
\item When $c<1$,  the map  $g$ has an unique maximum $s_1$ on $(0,\infty)$
  and an unique minimum $s_2$ on $(-\frac1{bc},0)$. 
  The edges of the support interval (for the absolutely continuous component)
  are again given by these local extrema:
  $x_1=g(s_1)$ and $x_2=g(s_2)$.
\end{enumerate}
\end{proposition}

\subsection{Application to an ARMA process}

For simplicity, we consider the simplest causal ARMA(1,1) process for the
coordinates:
\[    z_t=\phi z_{t-1}+  \veps_t+   \te \veps_{t-1},\quad t\in\Z,
\]
where $|\phi|<1$ and $\te$ is real.
The aim is to find a simplified form of general equation~\eqref{sol}. We
have 
\[     \frac1{2\pi f(\la)}= \left| \frac{1-\phi e^{i\la}}{1+\te e^{i\la}}\right|^2~,
\]
and 
\begin{eqnarray*}
  I &=& \frac{1}{2\pi}\int_{0}^{2\pi} \frac1{cs+\{2\pi f(\la)\}^{-1}}\,d\la
  =  \frac{1}{2\pi i}\oint_{|\xi|=1} 
      \frac1{cs+  \left| \frac{1-\phi \xi}{1+\te \xi}\right|^2  }   \frac{d\xi}{\xi}.
\end{eqnarray*}
By a lengthy but elementary calculation of residues detailed in
Section~\ref{proofs}, we find 
\begin{equation}\label{I}
  I = \frac{\te}{cs\te-\phi} 
  \left\{   1 - \frac{(\phi+\te)(1+\phi\te)}{\te(cs\te-\phi)}
  \frac{  \epsilon(\al)}{ \sqrt{\al^2-4}}
  \right\},
\end{equation}
with 
\begin{equation}\label{alpha}
  \al=\frac{cs(1+\te^2)+1+\phi^2}{cs\te-\phi} ,\qquad
  \epsilon(\al)=\sgn (\Im \al ) ~.
\end{equation}
Therefore for an ARMA(1,1) process, the general equation \eqref{sol} reduces  to 
\begin{equation}\label{sol-arma}
   z = -\frac1{s}  + \frac{\te}{cs\te-\phi} -   
   \frac{(\phi+\te)(1+\phi\te)}{(cs\te-\phi)^2}
   \frac{\epsilon(\al)}{\sqrt{\al^2-4}}~.
\end{equation}

Let us mention that it is important to have an explicit formula for  the integral 
in \eqref{sol}  to implement   numerical algorithms like the one
proposed in Section~\ref{sec:algo} in order to compute the 
density function of the LSD $F$.

\bigskip\noindent{\bf Case of an AR(1).}\quad For this particular
case, we have $\te=0$ and $\al=-(cs+1+\phi^2)/\phi$.
As $\Im s>0$, $\epsilon(\al)=\sgn(-\phi)$. It follows that 
\begin{eqnarray*}
  && - \frac{(\phi+\te)(1+\phi\te)}{(cs\te-\phi)^2}
  \frac{\epsilon(\al)}{\sqrt{\al^2-4}} 
  =  \frac{\sgn(-\phi)}{-\phi} \frac1{\sqrt{\frac{(cs+1+\phi^2)^2}{\phi^2}-4 }}\\
  &&\quad =  \frac{1}{|\phi|} \frac1{\sqrt{\frac{(cs+1+\phi^2)^2}{\phi^2}-4  }}
    =  \frac1{\sqrt{ (cs+1+\phi^2)^2 -4\phi^2  }}~.
\end{eqnarray*}
Therefore the Stieltjes
transform $s$ of the LSD is solution to a simpler equation
\begin{equation}\label{sol-ar}
   z = -\frac1{s}  +  \frac1{\sqrt{  (cs+1+\phi^2)^2-4\phi^2 }}~.
\end{equation}
It is worth noticing that if we further assume $\phi=0$, this  equation
reduces to $z=-1/s + 1/(cs+1)$ which characterises the standard \MP\ law
with i.i.d.  coordinates. 
Furthermore, for the determination of the support 
$[x_1,x_2]$ of the LSD, we notice that
\[  2\pi f(\la) = \frac1{|1-\phi e^{i\la}|^2}, 
\] 
so that its extrema are $a=1/(1+|\phi|)^2$ and 
$b=1/(1-|\phi|)^2$ (see Figure~\ref{fig:map-g}).

\bigskip\noindent{\bf Case of an MA(1).}\quad Here we have 
$\phi=0$ and 
\[  \al=\frac1\te\left(\frac1{cs}+1+\te^2\right)~.
\]
Hence $\epsilon(\al)=-\sgn(\te)$ and 
it  is readily checked out that the Stieltjes transform of the LSD 
is solution to the equation 
\begin{equation}\label{sol-ma}
  z = -\frac1s+\frac{1}{cs}  +  \frac1{c^2s^2} 
  \frac1{\sqrt{\left( \frac1{cs}+ 1+\te^2\right)^2 -4\te^2}}~.
\end{equation}
Again if we  further assume $\te=0$, this equation reduces to the one 
for the standard \MP\ law. 
Furthermore, for the determination of the support 
$[x_1,x_2]$ of the LSD, we notice that
\[  2\pi f(\la) = {|1+\te e^{i\la}|^2}, 
\] 
so that its extrema are  $a=(1-|\te|)^2$ and 
$b=(1+|\te|)^2$ (see Figure~\ref{fig:map-g}).

\section{A numerical method for computing the LSD  density function }
\label{sec:algo}
In this section we provide a numerical algorithm for the computation
of the density function $h$ of the LSD defined in Eq.\eqref{sol}
through its Stieltjes transform $s$.
We have
\[  s = \frac{1} {-z + A(s(z))}
\]
with 
\[    A(s(z))= \frac{1}{2\pi}\int_{0}^{2\pi}  \frac1{cs+\{2\pi f(\la)\}^{-1}}\,d\la~.
\]
The algorithm we propose is of fixed-point type. 
\begin{description}
\item[Algorithm]  
  For a given real $x$, let $\veps$ be small enough
  positive value and set $z=x+i\veps.$
 
  Choose an initial value $s_0(z)=u+i\veps$ and 
  iterate for $k\ge 0$ the above mapping  
  \[ s_{k+1}(z)=\{  -z+A(s_k(z))\}^{-1} ~,\]
  until convergence and let $s_K(z)$ be the final value.

  Define the estimate of the density function $h(x)$ to be 
  \[       \widehat h(x) = \frac1{\pi} \Im s_K(z)~.   \hskip2cm \blacksquare
  \]
\end{description}

It is well-known that  this iterated map has good
contraction properties  guaranteeing the convergence of the algorithm.
There are however two 
issues which need a careful consideration. First the integral 
operator $A$ is usually approximated by a numeric routine and because
of a high  number of calls to $A$, the resulting algorithm is slow. 
In this aspect analytic formula for $A$ when available are  well
acknowledged as  Eq.\eqref{sol-arma} in the case of  an  ARMA(1,1).

A second issue is that overall we first need to determine the support
interval $[x_1,x_2]$ of the density function $h$. This is handled with
the help of description of $x_j$'s given in 
Proposition~\ref{prop}. 

\begin{figure}[hb]
  \begin{center}
    \begin{tabular}{cc}
    \includegraphics[height=15pc,angle=-90]{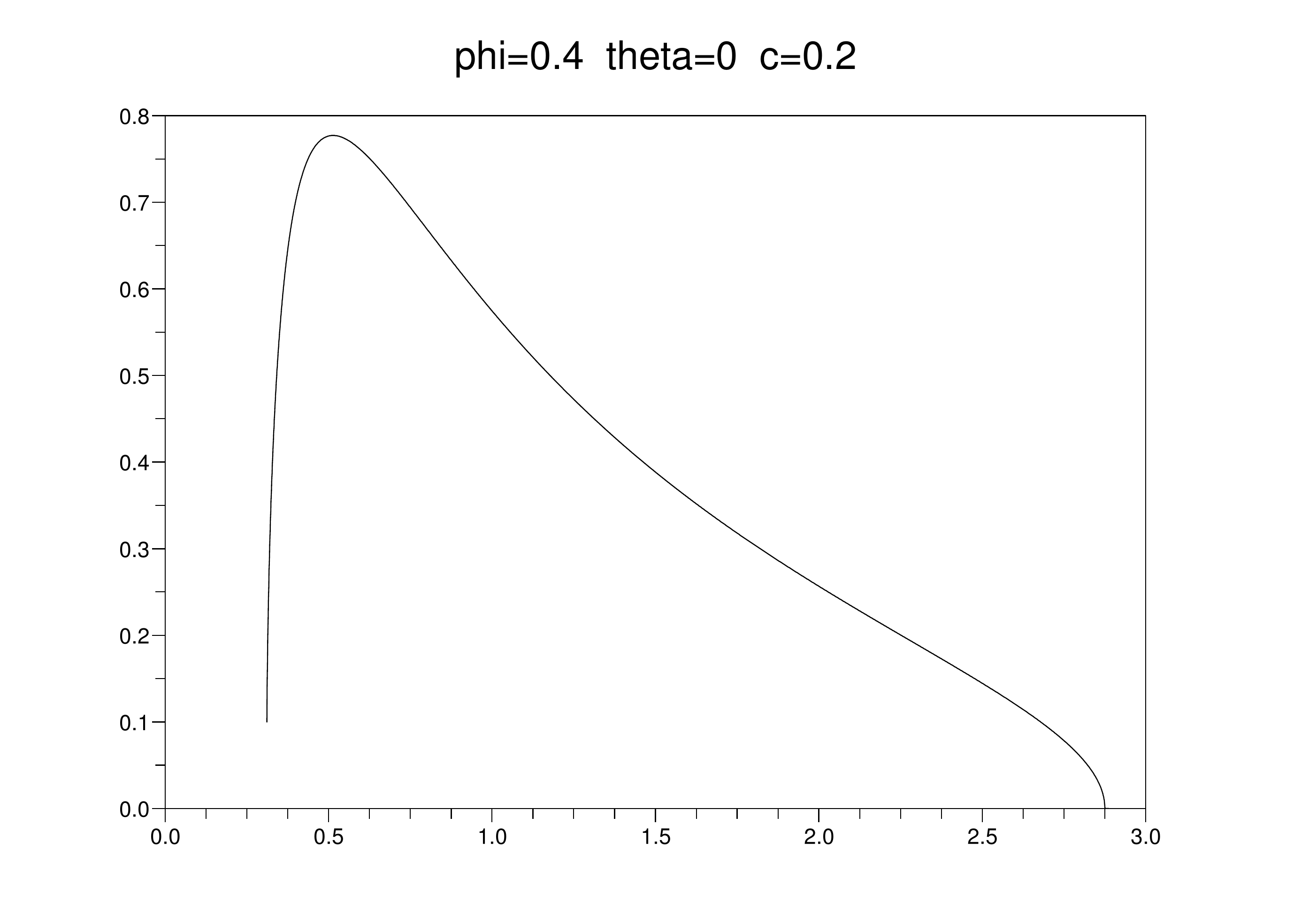}
&    \includegraphics[height=15pc,angle=-90]{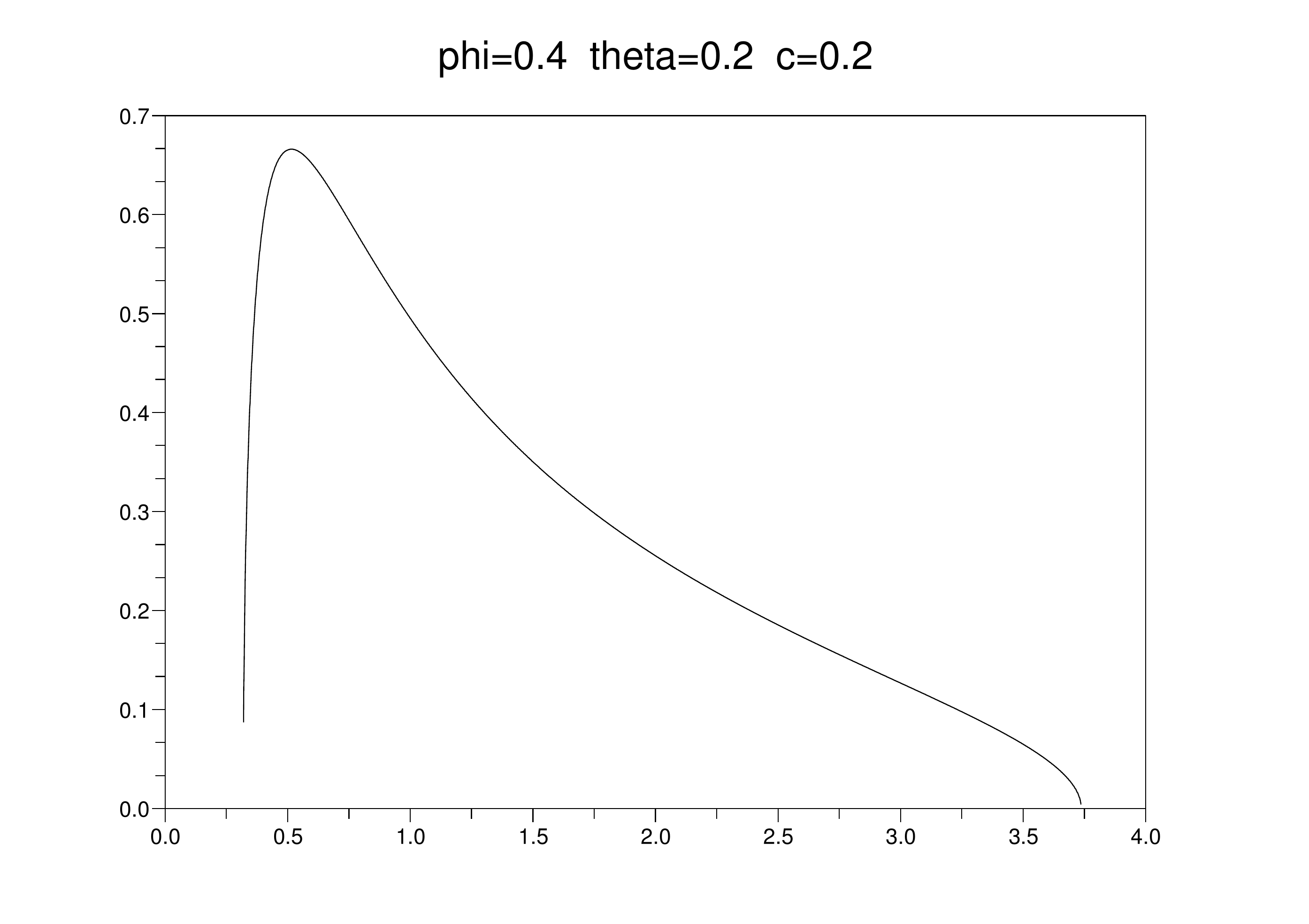}\\
    \includegraphics[height=15pc,angle=-90]{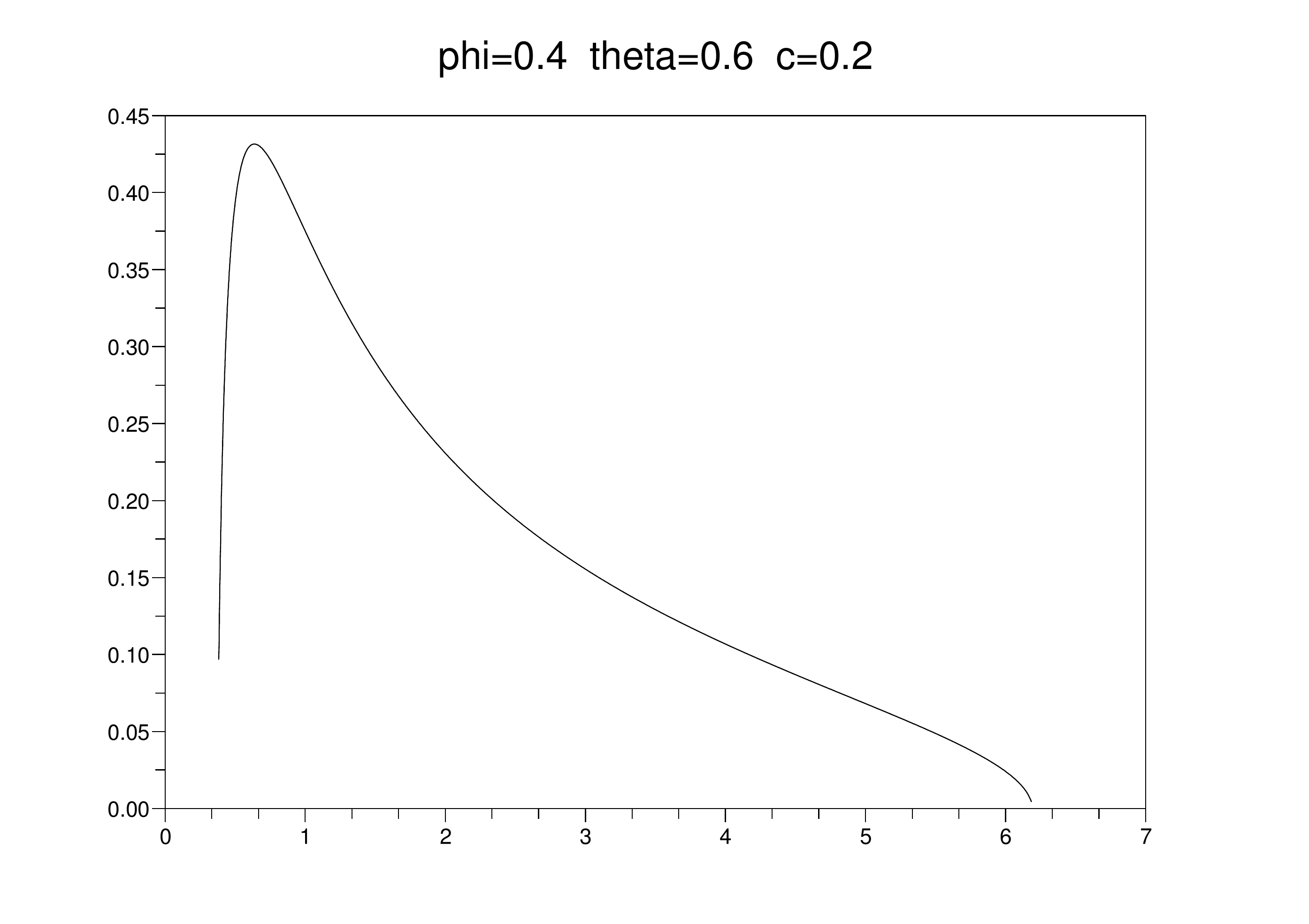}
&    \includegraphics[height=15pc,angle=-90]{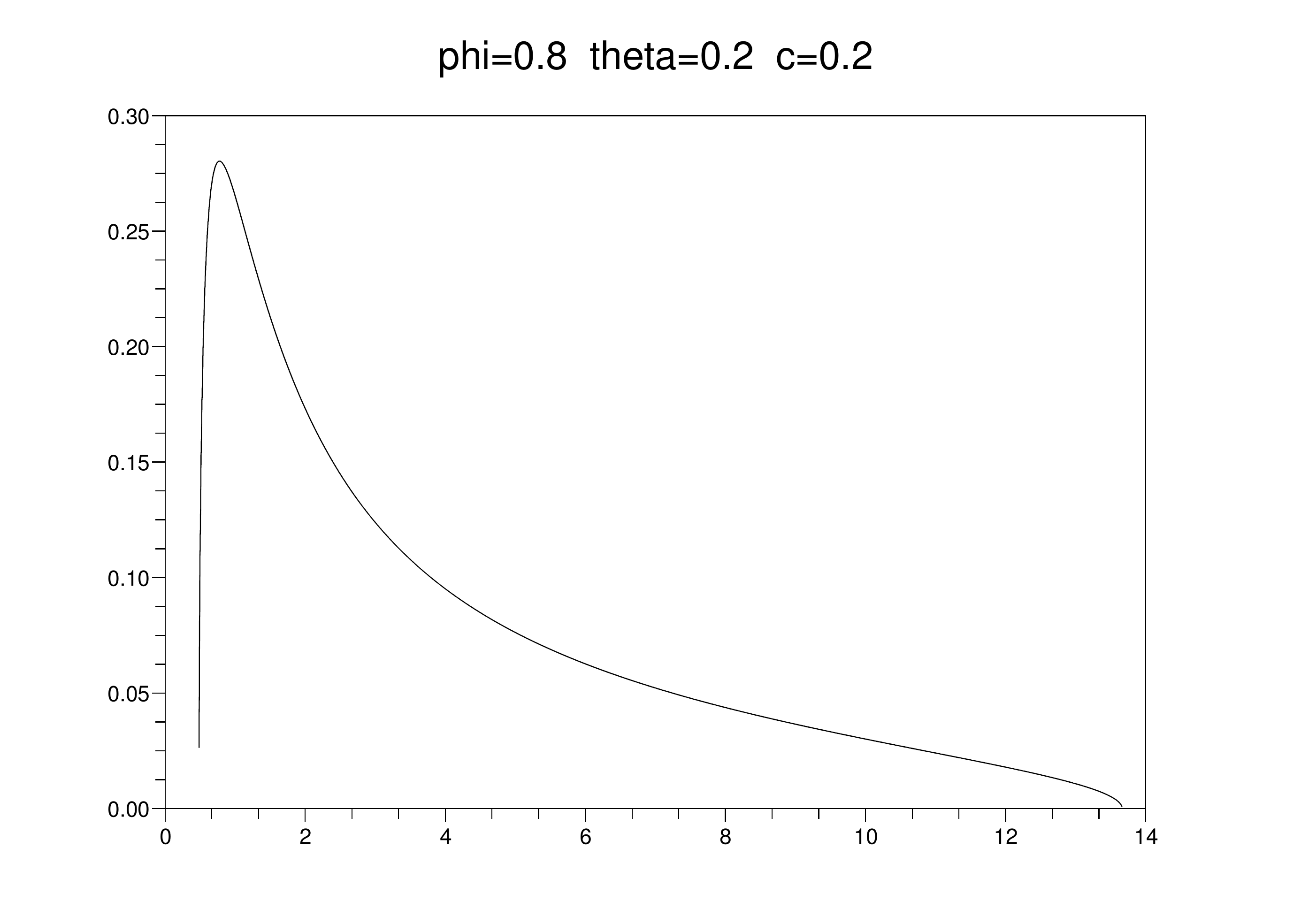}
    \end{tabular}
  \end{center}
  \caption{\label{fig:arma}Densities of the LSD from  ARMA(1,1)
    model. Left to right and top to bottom: $(\phi,\te,c)=$
    (0.4,~0,~0.2), 
    (0.4,~0.2,~0.2),  
    (0.4,~0.6,~0.2), and   
    (0.8,~0.2,~0.2).}  
\end{figure}

For four ARMA(1,1) models listed in
Table~\ref{tab:arma}, we have used this algorithm 
with the map $A$ defined in \eqref{sol-arma}
to get the density plots displayed in 
Figure~\ref{fig:arma}.

\renewcommand{\arraystretch}{1.5}
\begin{center}
  \begin{table}[h]
    \caption{\label{tab:arma}ARMA(1,1) models for density
      plots.  Reference support interval for  \MP\ law of index $c=0.2$
      is $[0.306, 2.094]$.}
    \begin{tabular}{ccc}
      Parameters $(\phi,\te,c)$   &  Estimated support  $[x_1,x_2]$\\
      \hline
      (0.4,~0,~0.2)  & [0.310,~2.875] \\
      (0.4,~0.2,~0.2)  & [0.319,~3.737] \\
      (0.4,~0.6,~0.2)  & [0.382,~6.186] \\
      (0.8,~0.2,~0.2)  & [0.485,~13.66] \\ \hline
    \end{tabular}
  \end{table}
\end{center}


Compared to the reference   standard
\MP\ law  with the same dimension to sample ratio $c$, 
 the above  density functions from ARMA models share a similar
 shape 
 with however a support interval getting larger and larger  with
 increasing 
 ARMA coefficients $\phi$ and $\te$.

\section{Proofs}\label{proofs}

\subsection*{Proof of Theorem~\protect\ref{th:main}}

Recall that the $p$ coordinates of the vectors 
$\bX_1,\ldots,\bX_n$ are i.i.d. while the temporal covariances
$\cov(X_{is},X_{it})$  are by definition those of $(z_t)$: for all
$1\le i\le p$, 
\[  \cov(X_{is},X_{it}) = \cov(z_s,z_t)=\ga_{t-s}, \quad 1\le s,t\le n~.
\]
Let $\widetilde f=2\pi f$.
It follows that 
the covariance matrix $T_n$ of each coordinate process 
$(X_{i1},\ldots,X_{in})$  equals to the   $n$-th order Toeplitz
matrix associated to $\widetilde f$:  
\[      T_n(s,t)= \ga_{t-s},\quad  \quad 1\le s,t\le n~,
\]
and 
\[    \widetilde f(\la) =\sum_{k=-\infty}^\infty \ga_k e^{ik\la}~,\quad \la\in[0,{2\pi}).      
\]   

We are going to apply Theorem~1.1 of \citet{BaiZhou08} for  
a strong limit of  the ESD of the sample covariance matrix 
$S_n=\frac1n \sum\limits_{t=1}^n \bX_t\bX_t^\top$.
Under the assumptions made, all the conditions of this theorem are satisfied 
except that we need to ensure a weak limit for   spectral
distributions of  $(T_n)$.  

The function $\widetilde f$ belongs to the Wiener class, i.e. the sequence
of its Fourier coefficients is absolutely summable. Moreover note that 
$\widetilde f$ is infinitely differentiable, its minimum $a$ and
maximum $b$ are attained. 
According to  the fundamental eigenvalue distribution theorem of
Szeg\"o for Toeplitz forms, see \citet[sect. 5.2]{GreSze58},
for any function $\varphi$ continuous on   
$[a, b]$ and denoting the eigenvalues of $T_n$ by
$\si_1^{(n)},\ldots,\si_n^{(n)}$, it holds that 
\[  \lim_{n\to\infty}\frac1n \sum_{k=1}^n \varphi(\si_k^{(n)}) =  
\frac1{2\pi} \int_{0}^{{2\pi}} \varphi(\widetilde f(\la))d\la~.
\]
Consequently,  the ESD of $T_n$ (i.e. distribution generated by the 
$\si_n^{(k)}$'s)  weakly converges to a nonrandom distribution 
$H$ with support $[a,b]$ and defined by 
\begin{equation}\label{H}
  H(x) = \frac1{2\pi} \int_{0}^{{2\pi}}  \mbox{\large 1}_{\{\widetilde
    f(\la)\le x\}} d\la~,
\end{equation}
and we have for  $\varphi$ as above, 
\begin{equation}\label{limitphi}
  \int_0^\infty \varphi(x) dH(x) = \frac1{2\pi} \int_{0}^{{2\pi}} \varphi(\widetilde f(\la))d\la~.
\end{equation}

Furthermore, by application of Theorem~1.1 of
\citet{BaiZhou08},  it holds that 
the ESD of $\frac{n}{p} S_n$ converges almost surely to a
nonrandom probability distribution whose  Stieltjes transform $m$
solves the equation 
\begin{eqnarray*}
  z & =&  - \frac1m + \frac1c  \int \frac{x}{1+mx} dH(x) \\
  &=&     - \frac1m + \frac1{2\pi c}   \int_{0}^{2\pi} \frac{1}{m+1/{\widetilde f}} d\la~, \\
\end{eqnarray*}
where we have used \eqref{limitphi} in the last equation.
The equation\eqref{sol} 
follows by observing the relation $s(z)=\frac1c m(z/c)$. 

\subsection*{Proof of Equation~\protect\eqref{I}}

The aim is to evaluate  the integral 
\begin{eqnarray*}
  &&I  =  \frac{1}{2\pi i}\oint_{|\xi|=1} 
  \frac1{cs+  \left| \frac{1-\phi \xi}{1+\te \xi}\right|^2 }  \frac{d\xi}{\xi}.
\end{eqnarray*}
Let $m=cs$ in this computation of residues.
We have
\begin{eqnarray*}
  I &=& \frac{1}{2\pi i}\oint_{|\xi|=1}     
  \frac{1+\te^2+\te(\xi+\xi^{-1})}
       { m\left\{1+\te^2+\te(\xi+\xi^{-1})\right\} + 
         \left\{ 1+\phi^2 -\phi(\xi+\xi^{-1} )\right\} } \frac{d\xi}{\xi}  \\
  &=&  \frac{\te}{m\te-\phi}   \frac{1}{2\pi i}\oint_{|\xi|=1}
       \left\{ \frac1\xi -  \frac{(\phi+\te)(1+\phi\te)}{\te(m\te-\phi)}
       \frac1{\xi^2+1+ \al \xi}
       \right\}
       d\xi ~\\
  &=&  \frac{\te}{m\te-\phi} \left\{ 1 - 
       \frac{(\phi+\te)(1+\phi\te)}{\te(m\te-\phi)} 
       \frac{1}{2\pi i}\oint_{|\xi|=1}  \frac{d\xi}{P(\xi)}  
      \right\}, 
\end{eqnarray*}
with  
\[  P(\xi)=\xi^2+ \al \xi+1,\quad
 \al=\frac{m(1+\te^2)+1+\phi^2}{m\te-\phi}.
\]
Let $\xi_1,\xi_2$ be the roots of $P(\xi)=\xi^2+1+ \al \xi$. Then 
\[ \frac1{P(\xi)} = \left( \frac1{\xi-\xi_1}-\frac1{\xi-\xi_2}   \right)\frac1{\xi_1-\xi_2}~.
\]
As $\xi_1\xi_2=1$, only one of the two poles is inside the unit
circle.   It is readily checked that 
if $\Im \al >0$, then $ |\xi_1|<|\xi_2|$ and 
\[\frac{1}{2\pi i}\oint_{|\xi|=1}  \frac{d\xi}{P(\xi)} = \frac1{\xi_1-\xi_2}=\frac1{\sqrt{\al^2-4}}.
\]
Otherwise we have $ |\xi_1|>|\xi_2|$ and the integral has an opposite
sign.  Summarising both cases we get 
\[  \frac{1}{2\pi i}\oint_{|\xi|=1}  \frac{d\xi}{P(\xi)} = \frac{\epsilon(\al)}{\sqrt{\al^2-4}},
\]
with $\epsilon(\al)=\sgn(\Im(\al))$. Equation~\eqref{I} is proved.

\bigskip
\noindent{\bf Acknowledgement.}   \quad 
The author is grateful to Jack Silverstein for several insightful discussions 
on the problem studied here, particularly for pointing to me  the numerical method 
of Section~\ref{sec:algo}.  We also thank a referee for important
comments on the paper.

\bibliographystyle{elsarticle-harv}

\vfill

\end{document}